\def\titlep{Continued fraction expansions and permutative representations 
of the Cuntz algebra ${\cal O}_{\infty}$}
\documentclass[11pt]{article}
\usepackage{graphicx,ifthen}
\usepackage{amssymb}

 at12pt

\newcommand{\qed}{\hbox{\rule[-2pt]{3pt}{6pt}}}
\newcommand{\qedh}{\hfill\qed \\}

\newtheorem{Thm}{Theorem}[section]

\newtheorem{rem}[Thm]{Remark}

\newtheorem{ex}[Thm]{Example}
\newtheorem{defi}[Thm]{Definition}

\newtheorem{prop}[Thm]{Proposition}

\newtheorem{fact}[Thm]{Fact}
\newtheorem{prob}[Thm]{Problem}

\newcommand{\vv}{\vspace{.3in}}
\newcommand{\ww}{\vv\noindent}

\def\cal#1{\mathcal #1}
\def\con{{\cal O}_{N}}
\def\coni{{\cal O}_{\infty}}

\def\pr{{\it Proof.}\quad}

\def\disp#1{{\displaystyle #1}}

\addtocounter{footnote}{1}
\def\sftt#1{
\setcounter{equation}{0}
\addtocounter{footnote}{1}
\section{#1}
}

\def\ssft#1{\subsection{#1}}

\def\sssft#1{\subsubsection{#1}}

\def\cls{\quad \clearpage}

\begin{document}
\def\inputcls#1{\cls\input #1.txt}
\def\plan#1#2{\par\noindent\makebox[.5in][c]{#1}
\makebox[.1in][l]{$|$}
\makebox[3in][l]{#2}\\}
\def\nset#1{\{1,\ldots,N\}^{#1}}
\def\brl{branching law}
\def\bfsnl{{\rm BFS}_{N}(\Lambda)}
\def\scm#1{S({\bf C}^{N})^{\otimes #1}}
\def\mqb{\{(M_{i},q_{i},B_{i})\}_{i=1}^{N}}
\newcommand{\mline}{\noindent
\thicklines
\setlength{\unitlength}{.18mm}
\begin{picture}(1000,5)
\put(0,0){\line(1,0){750}}
\end{picture}
}
\def\boxtimes{\noindent
\setlength{\unitlength}{.020918mm}
\begin{picture}(120,150)(60,0)
\thinlines
\put(0,0){
\line(1,0){100}\line(0,1){100}
}
\put(100,100){
\line(-1,0){100}\line(0,-1){100}
}
\put(40,0){$\times$}
\end{picture}
 }
\def\authorp{Katsunori  Kawamura}
\def\authorq{Dan Lascu}
\def\authorr{Ion Coltescu}
\def\authort{Yoshiki Hayashi}
\def\emailp{e-mail: kawamura@kurims.kyoto-u.ac.jp.}
\def\emailq{e-mail: lascudan@gmail.com.}
\def\emailr{e-mail: icoltescu@yahoo.com.}
\def\emailt{e-mail: yoshiki.h@hy5.ecs.kyoto-u.ac.jp.}
\def\addressp{{\small {\it College of Science 
and Engineering Ritsumeikan University,}}\\
{\small {\it 1-1-1 Noji Higashi, Kusatsu, Shiga 525-8577, Japan,}}
}
\def\addresst{{\small {\it College of Science 
and Engineering Ritsumeikan University,}}\\
{\small {\it 1-1-1 Noji Higashi, Kusatsu, Shiga 525-8577, Japan}}
}
\def\addressq{{\small {\it 
Mircea cel Batran Naval Academy, 1 Fulgerului, 900218 Constanta, Romania}}}
\def\ba{\mbox{\boldmath$a$}}
\def\bb{\mbox{\boldmath$b$}}
\def\bc{\mbox{\boldmath$c$}}
\def\be{\mbox{\boldmath$e$}}
\def\bp{\mbox{\boldmath$p$}}
\def\bq{\mbox{\boldmath$q$}}
\def\bu{\mbox{\boldmath$u$}}
\def\bv{\mbox{\boldmath$v$}}
\def\bw{\mbox{\boldmath$w$}}
\def\bx{\mbox{\boldmath$x$}}
\def\by{\mbox{\boldmath$y$}}
\def\bz{\mbox{\boldmath$z$}}
\def\aei{almost everywhere in}
\def\cdm#1{{\cal M}_{#1}(\{0,1\})}
\def\ptimes{\otimes_{\varphi}}
\def\cfe{{\sf CFE}}
\def\gpv{{\rm v}}
%
%
\pagestyle{plain}
\setcounter{page}{1}
\setcounter{section}{0}

\title{\titlep}

\author{
Katsunori Kawamura\thanks{\emailp},\,
\authort\thanks{\emailt}\\
\addressp\\
and \\
Dan Lascu\thanks{\emailq}\\
\addressq\\
}
\date{}

\maketitle
%
%
\begin{abstract}
We show a correspondence between 
simple continued fraction expansions of irrational numbers
and irreducible permutative representations 
of the Cuntz algebra ${\cal O}_{\infty}$.
With respect to the correspondence,
it is shown that 
the equivalence of real numbers with respect to modular transformations
is equivalent to the unitary equivalence of representations.
Furthermore, we show that 
quadratic irrationals are related to 
irreducible permutative representations of ${\cal O}_{\infty}$ with a cycle.
\end{abstract}

\sftt{Introduction}
\label{section:first}
The purpose of this paper is to show a new relation between 
number theory and the representation theory of operator algebras. 
At the beginning,
we show our motivation. 
Explicit mathematical statements will be given after 
$\S$ \ref{subsection:firsttwo}.
The main theorems will be shown in $\S$ \ref{subsection:firstfour}.

%
%
\ssft{Motivation}
\label{subsection:firstone}
Continued fractions furnish important tools 
in number theory \cite{HW,JT,RS}, and 
continued fraction transformations 
induce typical dynamical systems \cite{Chan,CL,IK,Schweiger}.
It is well known that the dynamical system
of the simple continued fraction transformation on the
set of irrational numbers in the interval $[0,1]$
is conjugate with the one-sided full shift on the set ${\bf N}^{\infty}
\equiv \{(n_{i})_{i\geq 1}:n_{i}\in {\bf N}\mbox{ for all }i\}$ 
\cite{Kitchens}
by using  continued fraction expansions
where ${\bf N}\equiv \{1,2,3,\ldots\}$.

On the other hand,
such dynamical systems induce representations of Cuntz algebras 
and their relations were studied \cite{PFO01,CFR01}.
For example, the one-sided full shift on ${\bf N}^{\infty}$ induces
the shift representation of $\coni$,
which acts on the representation space $l_{2}({\bf N}^{\infty})$ \cite{BJ}.
The shift representation is decomposed into the direct sum of irreducibles
unique up to unitary equivalence and the decomposition is multiplicity free.
Every irreducible component in the decomposition 
is a permutative representation, 
and any irreducible permutative representation
appears in the decomposition.

From these facts,
we are interested in a relation between
continued fraction expansions of irrationals and such representations 
of $\coni$ by the intermediary of the one-sided full shift on ${\bf N}^{\infty}$. 
We roughly illustrate relations among theories as follows:

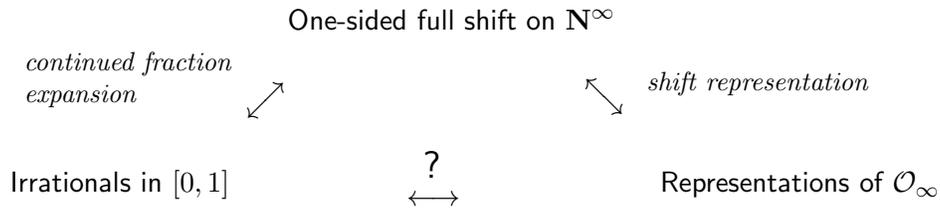
\def\figone{
{\sf
\put(20,0){
\begin{minipage}[t]{2in}
Irrationals in $[0,1]$
\end{minipage}
}
\put(225,120){
\begin{minipage}[t]{2in}
One-sided full shift on ${\bf N}^{\infty}$
\end{minipage}
}
\put(500,0){
\begin{minipage}[t]{2in}
Representations of $\coni$
\end{minipage}
}
\put(200,50){\rotatebox{45}{$\longleftrightarrow$}}
\put(30,90){
{\small {\it
\begin{minipage}[t]{2in}
continued fraction \\
expansion 
\end{minipage}
}
}
}
\put(450,60){\rotatebox{135}{$\longleftrightarrow$}}
\put(490,75){
{\small {\it
\begin{minipage}[t]{2in}
shift representation
\end{minipage}
}
}
}
\put(320,-10){$\longleftrightarrow$}
\put(332,10){{\Large ?}}
}
}
\noindent
\thicklines
\setlength{\unitlength}{.18mm}
\begin{picture}(800,180)(30,-30)
\put(0,0){\figone}
\end{picture}

\noindent
The position in the above question mark is the content
of this study.

%
%
\ssft{Continued fraction expansion map and continued fraction transformation
on the set of irrationals}
\label{subsection:firsttwo}
We review the continued fraction expansion map 
and the continued fraction transformation
according to \cite{HW}.
Let $[0,1]$ denote the closed interval from $0$ to $1$,
and let $\Omega$ denote the set of all irrationals in $[0,1]$, that is,
%
%
\begin{equation}
\label{eqn:omegaone}
\Omega= [0,1]\setminus {\bf Q}.
\end{equation}
Remark that we consider only $\Omega$ but not the whole of $[0,1]$
in this paper.
Any $x\in\Omega$ has a unique infinite continued fraction expansion
(\cite{HW}, Theorem 170),
that is, there exists a unique infinite sequence $(a_{i}(x))_{i\geq 1}$ 
of positive integers such that
%
%
\begin{equation}
\label{eqn:cfe}
x=\frac{1}{\disp{a_{1}(x)+\frac{1}{\disp{a_{2}(x)
+\frac{1}{a_{3}(x)+\disp{\frac{1}{\ddots}}}}}}}.
\end{equation}
From this, we define the map 
$\cfe$ from $\Omega$ to ${\bf N}^{\infty}$ by
%
%
\begin{equation}
\label{eqn:ax}
\cfe(x)\equiv (a_{i}(x))_{i\geq 1}.
\end{equation}
It is known that the map $\cfe$ is bijective 
(\cite{HW}, Theorem 161, 166, 169).
We call $\cfe$ {\it the continued fraction expansion map}.
If $x\in \Omega$ is a solution of a quadratic equation 
with integral coefficients,
then we call $x$ a {\it quadratic irrational}.
%
%
\begin{fact}
\label{fact:period}
(\cite{HW}, Theorem 177)
If $x\in \Omega$ is a quadratic irrational,
then there exists a finite sequence $(m_{1},\ldots,m_{l})\in {\bf N}^{l}
\cup\{\emptyset\}$ and 
a nonperiodic finite sequence $(n_{1},\ldots,n_{k})\in {\bf N}^{k}$ such that
%
%
\begin{equation}
\label{eqn:periodic}
\cfe(x)=(m_{1},\ldots,m_{l},n_{1},\ldots,n_{k},n_{1},\ldots,n_{k},
n_{1},\ldots,n_{k},\ldots).
\end{equation}
\end{fact}

\noindent
In Fact \ref{fact:period},
we call $(n_{1},\ldots,n_{k})$ the {\it repeating block} 
of $\cfe(x)$ (\cite{RS}, Chap. 3)
where we suppose that there  is no shorter such repeating block 
and that the initial block does not end with a copy of the 
repeating block.
Hence the repeating block of $\cfe(x)$ is uniquely defined
for each quadratic irrational $x$ in $\Omega$.
The converse of Fact \ref{fact:period} is also true (\cite{HW}, Theorem 176).

%
%
\ssft{Permutative representations of $\coni$}
\label{subsection:firstthree}
We review permutative representations of the Cuntz algebra $\coni$ 
in this subsection.
%
%
\sssft{$\coni$}
\label{subsubsection:firstthreeone}
Let $\coni$ denote the {\it Cuntz algebra} \cite{C}, that is, 
a C$^{*}$-algebra which is universally generated 
by $\{s_{i}:i\in {\bf N}\}$ satisfying
%
%
\begin{equation}
\label{eqn:coni}
s_{i}^{*}s_{j}=\delta_{ij}I\quad(i,j\in {\bf N}),\quad 
\sum_{i=1}^{k}s_{i}s_{i}^{*}\leq I,\quad (k\in {\bf N}),
\end{equation}
where $I$ denotes the unit of $\coni$.

A {\it $*$-representation} of $\coni$ is a pair $({\cal H},\pi)$
such that 
$\pi$ is a $*$-homomorphism from $\coni$ 
to the C$^{*}$-algebra ${\cal L}({\cal H})$
of all bounded linear operators 
on a complex Hilbert space ${\cal H}$ \cite{Blackadar2006}.
We call $*$-representation as representation for the simplicity of description.
For two representations $({\cal H}_{1},\pi_{1})$ and 
$({\cal H}_{2},\pi_{2})$ of $\coni$,
$({\cal H}_{1},\pi_{1})$ and 
$({\cal H}_{2},\pi_{2})$ are {\it unitarily equivalent}
if there exists a unitary $u$ from ${\cal H}_{1}$ onto ${\cal H}_{2}$
such that $u\pi_{1}(x) u^{*}=\pi_{2}(x)$ for each $x\in \coni$.
We state that a representation $({\cal H},\pi)$ of $\coni$ 
is {\it irreducible} if 
there is no invariant closed subspace of ${\cal H}$
except $\{0\}$ and ${\cal H}$;
$({\cal H},\pi)$ is {\it multiplicity free} 
if any two subrepresentations of $({\cal H},\pi)$
are not unitarily equivalent.
A representation $({\cal H},\pi)$ is irreducible if and only if
the {\it commutant} $\{X\in {\cal L}({\cal H}):X\pi(A)=\pi(A)X
\mbox{ for all }A\in \coni\}$ of $\pi(\coni)$ equals to ${\bf C}I$.

Since $\coni$ is simple, that is, there is no
nontrivial closed two-sided ideal,
any representation of $\coni$ is injective.
If $\{t_{i}:i\in {\bf N}\}$ are bounded operators 
on a Hilbert space ${\cal H}$ such that
$\{t_{i}:i\in {\bf N}\}$ satisfy (\ref{eqn:coni}),
then the correspondence $s_{i}\mapsto t_{i}$ for $i\in {\bf N}$
is uniquely extended to a unital $*$-representation
of $\coni$ on ${\cal H}$ from the uniqueness of $\coni$.
Therefore we call such a correspondence 
among generators by a representation of $\coni$ on ${\cal H}$.
Assume that $\{s_{i}:i\in {\bf N}\}$ are realized as
operators on a Hilbert space ${\cal H}$.
According to (\ref{eqn:coni}),
${\cal H}$ is decomposed into orthogonal subspaces as
$\oplus_{i\in {\bf N}}s_{i}{\cal H}$.
Since $s_{i}$ is an isometry,
$s_{i}{\cal H}$ has the same dimension as ${\cal H}$.
From this, we see that 
there is no finite dimensional representation of $\coni$ which preserves 
the unit.
The following illustration is helpful in understanding 
$\{s_{i}:i\in {\bf N}\}$:

%
%
\def\firstbox{
\put(0,0){\line(1,0){500}}
\put(0,30){\line(1,0){500}}
\put(0,0){\line(0,1){30}}
\put(500,0){\line(0,1){30}}
\put(240,10){${\huge {\cal H}}$}
}
\def\secondbox{
\put(0,0){\line(1,0){500}}
\put(0,30){\line(1,0){500}}
\put(0,0){\line(0,1){30}}
\put(100,0){\line(0,1){30}}
\put(200,0){\line(0,1){30}}
\put(300,0){\line(0,1){30}}
\put(140,10){$\cdots$}
\put(340,10){$\cdots$}
\put(500,0){\line(0,1){30}}
\put(40,10){${\huge s_{1}{\cal H}}$}
\put(240,10){${\huge s_{i}{\cal H}}$}
}
\def\cross{
\thinlines
\put(225,60){$s_{i}$}
\put(245,60){$\downarrow$}
\qbezier[200](500,100)(400,65)(300,30)
\qbezier[200](0,100)(100,65)(200,30)
}
\setlength{\unitlength}{.22mm}
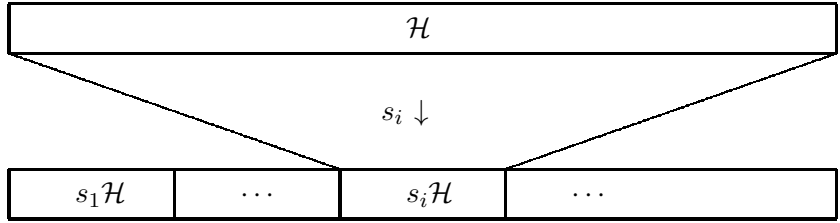
\begin{picture}(1000,150)(0,-10)
\thicklines
\put(0,100){\firstbox}
\put(0,0){\secondbox}
\put(0,0){\cross}
\end{picture}

The algebra $\coni$ appears in quantum field theory \cite{RBS01}
and metrical number theory \cite{CFR01}.

%
%
\sssft{Permutative representations}
\label{subsubsection:firstthreetwo}
We review permutative representations in this subsubsection.
%
%
\begin{defi}\cite{BJ,DaPi2,DaPi3,RBS01}
\label{defi:first}
Let $\{s_{i}:i\in {\bf N}\}$ denote the canonical generators of $\coni$.
\begin{enumerate}
\item
A representation $({\cal H},\pi)$ of $\coni$ is permutative if 
there exists an orthonormal basis ${\cal E}(\subset {\cal H})$ 
of ${\cal H}$
such that $\pi(s_{i}){\cal E}\subset {\cal E}$ for each $i\in {\bf N}$.
\item
For $J=(j_{l})_{l=1}^{k}\in {\bf N}^{k}$ with $1\leq k < \infty$,
let $P(J)$ denote 
the class of representations $({\cal H}, \pi)$ of $\coni$ 
with a cyclic unit vector $\gpv\in {\cal H}$
such that $\pi(s_{J})\gpv=\gpv$
and $\{\pi(s_{j_{l}}\cdots s_{j_{k}})\gpv\}_{l=1}^{k}$
is an orthonormal family in ${\cal H}$
where $s_{J}\equiv s_{j_{1}}\cdots s_{j_{k}}$.
\item
For $J=(j_{l})_{l\geq 1}\in {\bf N}^{\infty}$,
let $P(J)$ denote 
the class of representations $({\cal H}, \pi)$ of $\coni$ 
with a cyclic unit vector $\gpv\in {\cal H}$ such that
$\{\pi(s_{J_{(n)}})^{*}\gpv:n\in {\bf N}\}$
is an orthonormal family in ${\cal H}$ where
$J_{(n)}\equiv (j_{1},\ldots,j_{n})$. 
\end{enumerate}
The vector $\gpv$ in both (ii) and (iii) is 
called the GP vector of $({\cal H},\pi)$.
\end{defi}

We recall properties of these classes as follows:
For any $J$, $P(J)$ in Definition \ref{defi:first}(ii) and (iii)
always exists and it is a class of permutative representations, which 
contains only one unitary equivalence class.
From this, we can always identify $P(J)$ with a representative of $P(J)$
\cite{BJ,DaPi2,DaPi3}.
A representation 
$({\cal H},\pi)$ of $\coni$ is a {\it
permutative representation with a cycle (chain)} 
if there exists $J\in {\bf N}^{k}$ for $1\leq k<\infty$
({\it resp.} $J\in {\bf N}^{\infty}$)
such that $({\cal H},\pi)$ is $P(J)$.
Details will be explained in $\S$ \ref{subsection:secondtwo}.

%
%
\ssft{Main theorems}
\label{subsection:firstfour}
We show our main theorems in this subsection.
For this purpose,
we construct two representations of $\coni$ as follows.
For a nonempty set $A$,
let $l_{2}(A)$ denote the complex Hilbert space
with an orthonormal basis $\{e_{a}:a\in A\}$.
We call $\{e_{a}:a\in A\}$ the {\it standard basis} of $l_{2}(A)$.
%
%
\begin{defi}
\label{defi:representation}
Let $\Omega$ be as in (\ref{eqn:omegaone}) and 
let $\{s_{i}:i\in {\bf N}\}$ denote the canonical generators of $\coni$.
\begin{enumerate}
\item
For $i\in {\bf N}$,
define the map $\alpha_{i}$ from $\Omega$ to $\Omega$ by
%
%
\begin{equation}
\label{eqn:functiontwo}
\alpha_{i}(x)\equiv \frac{1}{x+i}\quad(x\in \Omega).
\end{equation}
Define the representation $\pi_{\alpha}$ of $\coni$ on $l_{2}(\Omega)$ by
%
%
\begin{equation}
\label{eqn:snn}
\pi_{\alpha}(s_{i})e_{x}\equiv 
e_{\alpha_{i}(x)}\quad (x\in \Omega,\,i\in {\bf N}).
\end{equation}
\item
For $i\in {\bf N}$, define the map
$\beta_{i}$ from ${\bf N}^{\infty}$ to ${\bf N}^{\infty}$ by
%
%
\begin{equation}
\label{eqn:gfunction}
\beta_{i}(n_{1},n_{2},\ldots)\equiv (i,n_{1},n_{2},\ldots)
\quad((n_{1},n_{2},\ldots)\in {\bf N}^{\infty},\,i\in {\bf N}).
\end{equation}
Define the representation $\pi_{\beta}$ of $\coni$ 
on $l_{2}({\bf N}^{\infty})$ by
%
%
\begin{equation}
\label{eqn:sn}
\pi_{\beta}(s_{i})e_{a}\equiv e_{\beta_{i}(a)}
\quad (a\in {\bf N}^{\infty},\,i\in {\bf N}).
\end{equation}
The representation 
$(l_{2}({\bf N}^{\infty}),\pi_{\beta})$
is called the {\it shift representation} \cite{BJ}.
\end{enumerate}
\end{defi}

\noindent
Then the following holds.
%
%
\begin{Thm}
\label{Thm:formal}
Two representations $(l_{2}(\Omega),\pi_{\alpha})$
and $(l_{2}({\bf N}^{\infty}),\pi_{\beta})$
are unitarily equivalent.
\end{Thm}

\noindent
From Theorem \ref{Thm:formal},
we can compare irreducible components of 
$(l_{2}(\Omega),\pi_{\alpha})$ with those 
of $(l_{2}({\bf N}^{\infty}),\pi_{\beta})$
and consider how an irreducible component 
in $(l_{2}(\Omega),\pi_{\alpha})$ is realized.
For this purpose,
we introduce an equivalence relation of real numbers as follows.
%
%
\begin{defi}
\label{defi:modular}\cite{Coppel,HW,RS}
If $x$ and $y$ are two real numbers such that
%
%
\begin{equation}
\label{eqn:modulardef}
x=\frac{ay+b}{cy+d}\quad(a,b,c,d\in {\bf Z})
\end{equation}
where $ad-bc=\pm 1$,
then $x$ is said to be equivalent to $y$.
In this case, we write $x\sim y$.
The transformation 
(\ref{eqn:modulardef}) is called a modular transformation
in a broad sense.
\end{defi}

\noindent
Remark that the transformation (\ref{eqn:modulardef}) 
for $a,b,c,d$ with $ad-bc=1$
is also called a modular transformation in a narrow sense.

According to the unitary equivalence in Theorem \ref{Thm:formal}, 
the following holds.
%
%
\begin{Thm}
\label{Thm:formaltwo}
For $x\in\Omega$,
let $[x]$ denote the equivalence class of $x$ in $\Omega$
with respect to modular transformations, that is,
$[x]=\{y\in \Omega:y\sim x\}$.
Let $(l_{2}(\Omega),\pi_{\alpha})$ be as in Definition \ref{defi:representation}(i).
\begin{enumerate}
\item
The following irreducible decomposition holds:
%
%
\begin{equation}
\label{eqn:tomega}
l_{2}(\Omega)=\bigoplus_{[x]\in \Omega/\!\sim}{\cal H}_{[x]}
\end{equation}
where ${\cal H}_{[x]}$ denotes the closed subspace of $l_{2}(\Omega)$
generated by the subset $\{e_{y}:y\in [x]\}$.
\item
For $x\in \Omega$,
let $\eta_{[x]}$ denote the subrepresentation
of $\pi_{\alpha}$ associated with the subspace ${\cal H}_{[x]}$
in (\ref{eqn:tomega}), that is,
%
%
\begin{equation}
\label{eqn:eta}
\eta_{[x]}\equiv \pi_{\alpha}|_{{\cal H}_{[x]}}.
\end{equation}
Then $\eta_{[x]}$ and $\eta_{[y]}$ are unitarily equivalent 
if and only if $x\sim y$.
Especially, (\ref{eqn:tomega}) is multiplicity free.
\item
Let $\cfe$ be as in (\ref{eqn:ax})
and let $\Omega^{(2)}$ denote
the set of all quadratic irrationals in $\Omega$.
\begin{enumerate}
\item
If $x\in\Omega\setminus \Omega^{(2)}$, then
$\eta_{[x]}$ is $P(\cfe(x))$.
\item
If $x\in\Omega^{(2)}$, then
$\eta_{[x]}$ is $P(\cfe_{0}(x))$
where $\cfe_{0}(x)$ denotes the repeating block of $\cfe(x)$.
\end{enumerate}
\item
Any irreducible permutative representation of $\coni$
is unitarily equivalent to 
$\eta_{[x]}$ for some $[x]\in\Omega/\!\!\sim$.
\end{enumerate}
\end{Thm}

In consequence,
Theorem \ref{Thm:formaltwo} shows that the set $\Omega/\!\!\!\sim$
of all equivalence classes of irrationals in $[0,1]$ is 
one-to-one correspondence in 
the set ${\rm IPR}(\coni)/\!\!\!\sim$
of all unitary equivalence classes of irreducible permutative representations
of $\coni$:
%
%
\begin{equation}
\label{eqn:correspondence}
\Omega/\!\!\sim\quad \stackrel{\eta}{\cong}\quad {\rm IPR}(\coni)/\!\!\sim;
\quad [x]\mapsto \eta_{[x]}
\end{equation}
where we identify $\eta_{[x]}$
with the unitary equivalence class of $\eta_{[x]}$ for convenience.
Especially, the following equivalence holds as the
restriction of $\eta$ on the subset $\Omega^{(2)}/\!\!\sim$
of $\Omega/\!\!\sim$:
%
%
\begin{equation}
\label{eqn:correspondencetwo}
\Omega^{(2)}/\!\!\sim\quad \stackrel{\eta}{\cong}\quad {\rm IPR}_{cycle}
(\coni)/\!\!\sim
\end{equation}
where ${\rm IPR}_{cycle}
(\coni)/\!\!\sim$ denotes the set of 
all unitary equivalence classes of irreducible permutative representations
of $\coni$ with a cycle.
%
%
\begin{rem}
\label{rem:one}
{\rm
The equivalence  of numbers by modular transformations
is well known in number theory, such that 
the discriminant of an irrational is invariant 
with respect to modular transformations 
\cite{Coppel,Zagier}.
On the other hand, the unitary equivalence of representations is
basic in the representation theory of $*$-algebras.
Since these two equivalence relations are independently
introduced in different mathematical areas, 
the equivalence of two equivalence relations
in Theorem \ref{Thm:formaltwo}(ii)
is nontrivial.
}
\end{rem}

From Theorem \ref{Thm:formaltwo},
the following natural questions are thought up.

%
%
\begin{prob}
\label{prob:one}
{\rm
\begin{enumerate}
\item
Show the meaning of the discriminant 
for the representation associated with a quadratic irrational.
What are the discriminant and the class number (\cite{Zagier}, p. 59)
in the representation theory of $\coni$?
\item
Find similar relations between
the Cuntz algebra $\con$ with $2\leq N<\infty$ \cite{C}
and real numbers.
\item
From Theorem \ref{Thm:formaltwo},
we suspect that there exists a relation between 
the group of all modular transformations and $\coni$.
Make clear this relation.
\end{enumerate}
}
\end{prob}

In $\S$ \ref{section:second}, we prove Theorem \ref{Thm:formal}
and \ref{Thm:formaltwo}.
In $\S$ \ref{section:third}, we show examples of Theorem \ref{Thm:formaltwo}.
%
%
\sftt{Proofs of theorems}
\label{section:second}
In this section, we prove main theorems.
%
%
\ssft{Continued fraction expansion and one-sided full shift}
\label{subsection:secondone}
We review relations between continued fraction expansions 
and the one-sided full shift on ${\bf N}^{\infty}$ in this subsection.
Let  $[a_{1}(x),a_{2}(x),\ldots]$ denote (\ref{eqn:cfe}) 
for simplicity of description.
%
%
\begin{fact}
\label{fact:hwtwo}
(\cite{HW}, Theorem 175)
For the equivalence in Definition \ref{defi:modular},
two irrational numbers $x$ and $y$ 
are equivalent if and only if
there exist positive integers
$a_{1},\ldots,a_{m}$,
$b_{1},\ldots,b_{n}$ such that
%
%
\begin{equation}
\label{eqn:sequence}
x=[a_{1},\ldots,a_{m},c_{1},c_{2},\ldots],\quad
y=[b_{1},\ldots,b_{n},c_{1},c_{2},\ldots].
\end{equation}
\end{fact}

Define the {\it simple continued fraction transformation 
(or the Gauss map)} $\tau$ from $\Omega$ to $\Omega$ by
%
%
\begin{equation}
\label{eqn:gauss}
\tau(x)\equiv \frac{1}{x}-\left\lfloor \frac{1}{x}\right\rfloor
\quad(x\in \Omega),
\end{equation}
where $\lfloor\cdot \rfloor$ denotes the floor (entire) function
\cite{IK,Schweiger}. 
Then we see that
$\tau([a_{1},a_{2},\ldots])=[a_{2},a_{3},\ldots]$
for the continued fraction $[a_{1},a_{2},\ldots]$. 

Define the map 
$\sigma$ from ${\bf N}^{\infty}$ to ${\bf N}^{\infty}$ by 
%
%
\begin{equation}
\label{eqn:sigma}
\sigma(n_{1},n_{2},\ldots)\equiv (n_{2},n_{3},\ldots).
\end{equation}
The dynamical system $({\bf N}^{\infty},\sigma)$
is called the {\it one-sided full shift} on ${\bf N}^{\infty}$
(\cite{Kitchens}, $\S$ 7.2).
Then we see that $\cfe$ in (\ref{eqn:ax}) satisfies 
the following equation (\cite{IK}, (1.1.2)):
%
%
\begin{equation}
\label{eqn:conjugate}
\cfe(\tau(x))=\sigma(\cfe(x))\quad(x\in \Omega).
\end{equation}
From (\ref{eqn:conjugate}),
two dynamical systems $(\Omega,\tau)$
and $({\bf N}^{\infty},\sigma)$ are conjugate.
%
%
\begin{defi}
\label{defi:tail}
For $a,b\in {\bf N}^{\infty}$,
let $a\sim b$ denote
when there exist $p,q\geq 1$
such that $\sigma^{p}(a)=\sigma^{q}(b)$.
\end{defi}
We call $\sim$ the {\it tail equivalence} in ${\bf N}^{\infty}$ 
(\cite{BJ}, Chap. 2).
Then the following holds by definition.
%
%
\begin{fact}
\label{fact:equivalence}
For $x,y\in \Omega$,
$x\sim y$ if and only if $\cfe(x)\sim \cfe(y)$. 
\end{fact}
For $\{\alpha_{i}:i\in {\bf N}\}$ and 
$\{\beta_{i}:i\in {\bf N}\}$ in Definition \ref{defi:representation},
the following holds:
%
%
\begin{equation}
\label{eqn:conjugateb}
\cfe(\alpha_{i}(x))=\beta_{i}(\cfe(x))\quad(x\in \Omega,\,i\in {\bf N}).
\end{equation}
Remark that they are closely related to $\tau$ and $\sigma$ as follows:
%
%
\begin{equation}
\label{eqn:relations}
(\tau\circ \alpha_{i})(x)=x,\quad 
(\sigma\circ \beta_{i})(a)=a\quad 
(x\in \Omega,\,a\in {\bf N}^{\infty},\,i\in {\bf N}).
\end{equation}

%
%
\ssft{Properties of representations of $\coni$}
\label{subsection:secondtwo}
In this subsection,
we recall properties of permutative representations 
and the shift representation of $\coni$.
Define ${\bf N}^{*}\equiv \bigcup_{1\leq k<\infty}{\bf N}^{k}$.
For $J=(j_{l})_{l=1}^{m},K=(k_{l})_{l=1}^{m^{'}}\in {\bf N}^{*}$,
we write $J\sim K$ 
if $m=m^{'}$ and $K=pJ$
where  $p J=(j_{p(1)},\ldots,j_{p(m)})$ for any
cyclic permutation $p\in {\bf Z}_{m}$.
For $J\in {\bf N}^{*}$,
we call $J$ {\it nonperiodic} if
$p J\ne J$ for any cyclic permutation $p\ne id$.
If $J\in {\bf N}^{*}$ is nonperiodic,
then we see that
$({\cal H},\pi)$ is $P(J)$
if and only if there  exists a cyclic vector $\gpv\in {\cal H}$
such that $\pi(s_{J})\gpv=\gpv$.
For $J\in {\bf N}^{\infty}$,
we call $J$ {\it nonperiodic} if $J$ has no repeating block.
%
%
\begin{prop}
\label{prop:oniprop}
Let $P(J)$ be as in Definition \ref{defi:first}.
\begin{enumerate}
\item
Any permutative representation of $\coni$ is
decomposed into the direct sum of cyclic permutative representations 
uniquely up to unitary equivalence.
\item
Any cyclic permutative representation is 
either one of the following two cases:
\begin{enumerate}
\item
$P(J)$
for $J\in {\bf N}^{*}$. 
\item
$P(J)$ for $J\in {\bf N}^{\infty}$.
\end{enumerate}
\item
For two representations $\pi_{1}$ and $\pi_{2}$ of $\coni$,
let $\pi_{1}\sim \pi_{2}$ denote 
the unitary equivalence between 
$\pi_{1}$ and $\pi_{2}$.
If $J\in {\bf N}^{*}$ and $K\in {\bf N}^{\infty}$,
then $P(J)\not \sim P(K)$.
\item
For $J,K\in {\bf N}^{*}\cup {\bf N}^{\infty}$,
then $P(J)\sim P(K)$ if and only if $J\sim K$
where we define $J\not \sim K$
when $J\in {\bf N}^{*}$ and $K\in {\bf N}^{\infty}$.
\item
For $J\in {\bf N}^{*}\cup {\bf N}^{\infty}$,
$P(J)$ is irreducible if and only if $J$ is nonperiodic.
\end{enumerate}
\end{prop}
%
%
\pr
See Appendix \ref{subsection:apponeone}. \qedh

Next, we show properties of the shift representation of $\coni$
(see also \cite{BJ}, Chap. 6).
%
%
\begin{prop}
\label{prop:shift}
Let $(l_{2}({\bf N}^{\infty}),\pi_{\beta})$ be as in 
Definition \ref{defi:representation}(ii).
Define $[a]\equiv \{b\in {\bf N}^{\infty}:b\sim a\}$
where $\sim$ is as in Definition \ref{defi:tail}.
\begin{enumerate}
\item
The following irreducible decomposition holds:
%
%
\begin{equation}
\label{eqn:tomegatwo}
l_{2}({\bf N}^{\infty})
=\bigoplus_{[a]\in {\bf N}^{\infty}/\!\sim}{\cal K}_{[a]}
\end{equation}
where 
${\cal K}_{[a]}$ denotes the closed subspace of 
$l_{2}({\bf N}^{\infty})$ generated by the set $\{e_{b}:b\in [a]\}$.
\item
For $a\in {\bf N}^{\infty}$,
let $\theta_{[a]}$ denote the subrepresentation
of $\pi_{\beta}$ associated with the subspace ${\cal K}_{[a]}$, that is,
%
%
\begin{equation}
\label{eqn:etatwo}
\theta_{[a]}\equiv \pi_{\beta}|_{{\cal K}_{[a]}}.
\end{equation}
Then $\theta_{[a]}$ and $\theta_{[b]}$ are unitarily equivalent 
if and only if $a\sim b$.
Especially, (\ref{eqn:tomegatwo}) is multiplicity free.
\item
\begin{enumerate}
\item
If $a\in{\bf N}^{\infty}$ has no repeating block, then
$\theta_{[a]}$ is $P(a)$.
\item
If $a\in{\bf N}^{\infty}$ has the repeating block $a^{'}$, then
$\theta_{[a]}$ is $P(a^{'})$.
\end{enumerate}
\item
Any irreducible permutative representation of $\coni$
is unitarily equivalent to 
$\theta_{[a]}$ for some $[a]\in{\bf N}^{\infty}/\!\!\sim$.
\end{enumerate}
\end{prop}
%
%
\pr
See Appendix \ref{subsection:apponetwo}. \qedh

%
%
\ssft{Proofs of Theorem \ref{Thm:formal}\, 
and \ref{Thm:formaltwo}}
\label{subsection:secondthree}
We prove Theorem \ref{Thm:formal} and \ref{Thm:formaltwo}
in this subsection.\\

\noindent
{\it Proof of Theorem \ref{Thm:formal}.}
Let $\cfe$ be as in (\ref{eqn:ax}).
Define the unitary $U$ from $l_{2}(\Omega)$ to 
$l_{2}({\bf N}^{\infty})$ by
%
%
\begin{equation}
\label{eqn:unitary}
Ue_{x}\equiv e_{\cfe(x)}^{'}\quad (x\in\Omega)
\end{equation}
where we write $\{e_{x}:x\in \Omega\}$ and 
$\{e_{a}^{'}:a\in {\bf N}^{\infty}\}$
as standard basis of 
$l_{2}(\Omega)$  and $l_{2}({\bf N}^{\infty})$, respectively.
From (\ref{eqn:conjugateb}),
we can verify that $U\pi_{\alpha}(s_{i})U^{*}=\pi_{\beta}(s_{i})$
for each $i\in {\bf N}$.
This implies the unitary equivalence between
two representations $(l_{2}(\Omega),\pi_{\alpha})$ and 
$(l_{2}({\bf N}^{\infty}),\pi_{\beta})$:
%
%
\begin{equation}
\label{eqn:implement}
U\pi_{\alpha}(A)U^{*}=\pi_{\beta}(A)\quad(A\in \coni).
\end{equation}
Hence the statement holds.
\qedh

\ww
{\it Proof of Theorem \ref{Thm:formaltwo}.}
From Fact \ref{fact:equivalence},
we see that 
%
%
\begin{equation}
\label{eqn:conjugatezero}
\cfe([x])=[\cfe(x)]\quad(x\in\Omega).
\end{equation}
By this and (\ref{eqn:unitary}),
%
%
\begin{equation}
\label{eqn:unitarytwo}
U{\cal H}_{[x]}={\cal K}_{[\cfe(x)]}\quad(x\in\Omega).
\end{equation}
From Proposition \ref{prop:shift}(i)
and Theorem \ref{Thm:formal}, the statement of (i) holds.

By (\ref{eqn:implement}) and (\ref{eqn:unitarytwo}), 
we obtain that 
%
%
\begin{equation}
\label{eqn:unitarythree}
U\eta_{[x]}(A)U^{*}=\theta_{[\cfe(x)]}(A)\quad(x\in\Omega,\,A\in\coni).
\end{equation}
From Proposition \ref{prop:shift}(ii), (iii), (iv) and (\ref{eqn:unitarythree}),
the statements of (ii), (iii) and (iv) hold, respectively.
\qedh

%
%
\sftt{Permutative representations of $\coni$
associated with quadratic irrationals}
\label{section:third}
We show examples of Theorem \ref{Thm:formaltwo}(iii)-(b) in this section.
%
%
\begin{ex}
\label{ex:zero}
{\rm 
For $k\in {\bf N}$, define $x\in \Omega$ by
%
%
\begin{equation}
\label{eqn:onep}
x=\frac{\sqrt{k^{2}+4}-k}{2}.
\end{equation}
Then $x=[k,k,k,\dots]$ as an infinite continued fraction.
Therefore $\eta_{[x]}$ is $P(k)$.
For example,
%
%
\begin{equation}
\label{eqn:golden}
\eta_{[\frac{\sqrt{5}-1}{2}]}=P(1),\quad
\eta_{[\sqrt{2}-1]}=P(2),\quad
\eta_{[\frac{\sqrt{13}-3}{2}]}=P(3).
\end{equation}
In \cite{RBS01}, we showed that 
the restriction of $P(1)$ on the algebra ${\cal B}$ of bosons 
is unitarily equivalent to the Fock representation of ${\cal B}$.
Hence the first equation in (\ref{eqn:golden})
shows a relation between the golden ratio $(\sqrt{5}-1)/2$
and the Fock representation of bosons.
}
\end{ex}
%
%
\begin{ex}
\label{ex:four}
{\rm
For $j,k\in {\bf N}$,
define $x\in\Omega$ by
%
%
\begin{equation}
\label{eqn:twop}
x=\frac{\sqrt{(jk)^{2}+4jk}-jk}{2j}.
\end{equation}
When $j\ne k$,
%
%
\begin{equation}
\label{eqn:twoq}
\eta_{[x]}=P(j,k).
\end{equation}
For example,
$\eta_{[\sqrt{3}-1]}=P(1,2)$.
If $j=k$, then 
(\ref{eqn:twop}) equals to (\ref{eqn:onep}).
Therefore $\eta_{[x]}$ is $P(k)$.
}
\end{ex}
%
%
\begin{ex}
\label{ex:triple}
{\rm
For $i,j,k\in {\bf N}$,
define $x\in\Omega$ by
%
%
\begin{eqnarray}
\label{eqn:threep}
x=&\disp{\frac{-(ijk+i+k-j)+ \sqrt{D}}{2(ij+1)}},\\
\nonumber
\\ 
\label{eqn:threeptwo}
D=&(ijk+i+j+k)^{2}+4.
\end{eqnarray}
If $(i,j,k)$ is nonperiodic,
then $\eta_{[x]}=P(i,j,k)$.
If $i=j=k$, then 
(\ref{eqn:threep}) equals to (\ref{eqn:onep}).
}
\end{ex}
%
%
\begin{ex}
\label{ex:qualtet}
{\rm
For $i,j,k,l\in {\bf N}$,
define $x\in\Omega$ by
%
%
\begin{eqnarray}
\label{eqn:fourp}
x=&\disp{\frac{-(ijkl+ij+kl+li-jk)+ \sqrt{D}}{2(ijk+i+k)}},\\
\nonumber\\
\label{eqn:fivep}
D=&(ijkl+ij+jk+kl+li)(ijkl+ij+jk+kl+li+4).
\end{eqnarray}
If $(i,j,k,l)$ is nonperiodic,
then $\eta_{[x]}=P(i,j,k,l)$.
}
\end{ex}

\noindent
Remark that the symbol $D$ in 
neither (\ref{eqn:threeptwo}) nor (\ref{eqn:fivep})
always means the discriminant of $x$.
For example,
when $(i,j,k)=(1,2,3)$,
$D$ in (\ref{eqn:threeptwo})
is $148$.
On the other hand, 
$x$ in (\ref{eqn:threep}) is $\frac{-5+\sqrt{37}}{2}$
with the discriminant $37$.
%
%
\begin{prob}
\label{prob:last}
{\rm
For a given $J\in {\bf N}^{n}$ for $n\geq 5$,
compute $x\in\Omega$ such that 
$\eta_{[x]}=P(J)$.
}
\end{prob}

%

\appendix
\section*{Appendix}

%
%
\sftt{Proofs of propositions}
\label{section:appone}
We prove Proposition \ref{prop:oniprop} and \ref{prop:shift} in this section.
%
%
\ssft{Proof of Proposition \ref{prop:oniprop}}
\label{subsection:apponeone}
Let $({\cal H},\pi)$ be a permutative representation of $\coni$.
By assumption,
we see that 
there exists a family $\{f_{i}:i\in {\bf N}\}$ of maps on a set $\Lambda$
and an orthonormal basis $\{e_{n}:n\in \Lambda\}$
such that 
%
%
\begin{equation}
\label{eqn:permutative}
\pi(s_{i})e_{n}=e_{f_{i}(n)}\quad(i\in {\bf N},\,n\in\Lambda).
\end{equation}
From this and (\ref{eqn:coni}),
we see that
$\#\Lambda=\infty$,
the map $f_{i}:\Lambda\to \Lambda$ is injective for each $i$,
$f_{i}(\Lambda)\cap f_{j}(\Lambda)=\emptyset$ when $i\ne j$
and $\Lambda=\bigcup_{i\in {\bf N}}f_{i}(\Lambda)$.
The family $\{f_{i}:i\in {\bf N}\}$ is called a {\it branching function system} 
(\cite{BJ}, Definition 2.1).
On the other hand,
if a branching function system is given,
then we can construct a permutative representation of $\coni$
as (\ref{eqn:permutative}).
Hence we write $\pi$ in (\ref{eqn:permutative}) as $\pi_{f}$.
For a given branching function system $\{f_{i}:i\in {\bf N}\}$,
define the {\it coding map} $F$ 
of $\{f_{i}:i\in {\bf N}\}$ by the map from $\Lambda$ to $\Lambda$
as $F(n)\equiv f_{i}^{-1}(n)$ when $n\in f_{i}(\Lambda)$.
\\
\\
{\it Proof of Proposition \ref{prop:oniprop}.}
(i)
Let $({\cal H},\pi)$ be a permutative representation of $\coni$.
Then there exists a branching function system
$f=\{f_{i}:i\in {\bf N}\}$ on a set $X$ such that
 $({\cal H},\pi)$ is unitarily equivalent to 
$(l_{2}(X),\pi_{f})$ by definition.
For $x,y\in X$,
define the equivalence relation $\sim $ as
$x\sim y$ if and only if there exist $i,j\in {\bf N}$
such that $F^{i}(x)=F^{j}(y)$
where $F$ denotes the coding map of $f$.
Let $\{X_{\lambda}:\lambda\in \Xi\}$
denote the set of all equivalence classes in $X$ with respect to $\sim$.
Then we see that $X$ is decomposed into the disjoint union of subsets 
$\{X_{\lambda}:\lambda\in\Xi\}$
such that $f_{i}(X_{\lambda})\subset X_{\lambda}$ for
each $i\in {\bf N}$ and $\lambda\in \Xi$.
Therefore $f$ is decomposed into
the direct sum of branching function systems
$f^{(\lambda)}\equiv \{f_{i}|_{X_{\lambda}}:i\in {\bf N}\}$
for $\lambda\in\Xi$.
By the definition of $X_{\lambda}$,
$(l_{2}(X_{\lambda}),\pi_{f^{(\lambda)}})$ is cyclic
and $\pi_{f^{(\lambda)}}=\pi_{f}|_{l_{2}(X_{\lambda})}$.
This implies that 
${\cal H}$ is decomposed into the direct sum of 
cyclic subspaces $\{l_{2}(X_{\lambda}):\lambda\in\Xi\}$.
Hence the statement of decomposition holds.

Assume that $g=\{g_{i}:i\in {\bf N}\}$ 
is another branching function system
on a set $Y$ associated with $\pi$.
Since 
dimensions of the representation space of both $\pi_{f}$ 
and $\pi_{g}$ are same, we can assume that $Y=X$.
Define the map
$\phi$ from $X$ to $X$ by
$\phi(x)\equiv(g_{i}\circ f_{i}^{-1})(x)$
when $x\in f_{i}(X)$.
Then $\phi$ is bijective and induces the same decomposition up to conjugacy.
Define the unitary $U$ on ${\cal H}$ by
$Ue_{x}\equiv e_{\phi(x)}$ for the standard basis $\{e_{x}:x\in X\}$.
Then we see that $U\pi_{f}(\cdot)U^{*}=\pi_{g}$.
Hence the statement of uniqueness holds.

\noindent
(ii)
Let $({\cal H},\pi)$ be a cyclic permutative representation of $\coni$.
From the proof of (i),
we can assume that 
there exists a branching function system $f=\{f_{i}:i\in {\bf N}\}$
on a set $X$ such that 
$({\cal H},\pi)$ is unitarily equivalent to $(l_{2}(X),\pi_{f})$
with a cyclic vector $e_{x_{0}}$ for a point $x_{0}\in X$.
For the coding map $F$ of $f$,
define
the sequence $x_{n}\equiv F^{n}(x_{0})$ for $n\geq 0$.
If there exist $m_{0},k\geq 1$ such that $x_{m+k}=x_{m}$
for each $m\geq m_{0}$,
then we see that $({\cal H},\pi)$ is the case of (a).
If not, then $({\cal H},\pi)$ is the case of (b).

\noindent
(iii)
If they are equivalent,
then there exists an action of $\coni$ on
a Hilbert space ${\cal H}$
with two cyclic unit vectors $\gpv$ and $\gpv^{'}$
such that 
$\gpv$ and $\gpv^{'}$ are
GP vectors of $P(J)$ and $P(K)$, respectively.
Then
$\langle\gpv|\gpv^{'}\rangle
=\langle\gpv|s_{J}^{*}\gpv^{'}\rangle
=\cdots 
=\langle\gpv|(s_{J}^{*})^{n}\gpv^{'}\rangle
\to 0$ ($n\to \infty)$
because of the assumption of $P(K)$.
From this,
$\langle\gpv|\gpv^{'}\rangle=0$.
This implies that 
$\langle s_{L}\gpv|\gpv^{'}\rangle=0$
for each $L\in {\bf N}^{*}$.
Therefore $\gpv^{'}=0$.
This is a contradiction.

\noindent
(iv)
From (iii),
it is sufficient to show the following two cases.

\noindent
(a)
Assume $J,K\in {\bf N}^{*}$.
We can assume that $\coni$ acts on two Hilbert spaces ${\cal H}$
and ${\cal H}^{'}$
with cyclic unit vectors $\gpv,\gpv^{'}$ 
which are GP vectors of $P(J)$ and $P(K)$, respectively.
If $J\sim K$,
then we can find $\gpv^{''}\in {\cal H}^{'}$ such that
$s_{J}\gpv^{''}=\gpv^{''}$
and $\{s_{j_{k}}\cdots s_{j_{n}}\gpv^{''}:k=1,\ldots,n\}$
is an orthonormal family when $J=(j_{1},\ldots,j_{n})$.
Define the unitary $U$ from ${\cal H}$ to ${\cal H}^{'}$
by $Us_{L}\gpv=s_{L}\gpv^{''}$
for each $L\in {\bf N}^{*}$.
Then $U$ gives the unitary equivalence of 
${\cal H}$ and ${\cal H}^{'}$.

Assume that $P(J)\sim P(K)$, but $J\not \sim K$.
Then we can assume that $\coni$ acts on a Hilbert space 
${\cal H}$ with two cyclic unit vectors $\gpv$ and 
$\gpv^{'}$ which are GP vectors of $P(J)$ and $P(K)$, respectively.
Then
$\langle\gpv|\gpv^{'}\rangle
=\langle\gpv|(s_{J}^{*})^{n}s_{K}^{m}\gpv^{'}\rangle$
for each $n,m\geq 1$.
Since $J\not \sim K$,
we see that $(s_{J}^{*})^{n}s_{K}^{m}=0$ for some $n,m$.
Therefore 
$\langle\gpv|\gpv^{'}\rangle=0$.
This implies that 
$\langle s_{L}\gpv|\gpv^{'}\rangle=0$
for each $L\in {\bf N}^{*}$.
Therefore $\gpv^{'}=0$.
This is a contradiction.

\noindent
(b)
Assume $J,K\in {\bf N}^{\infty}$.
Then $P(J)\sim P(K)$ if and only if $J\sim K$
from the analogy of the case (a).

\noindent
(v)
We consider the following two cases.

\noindent
(a)
Assume $J\in {\bf N}^{*}$.
Assume that $J$ is nonperiodic
and $\coni$ acts on a Hilbert space ${\cal H}$
with a cyclic unit vector $\gpv$
which is the GP vector of $P(J)$.
Let $w\in {\cal H}$.
Since ${\cal H}$ is generated by
the set $\{s_{K}\gpv:K\in {\bf N}^{*}\}$,
we can write 
$w=\sum_{L}a_{L}s_{L}\gpv$
such that the set $\{s_{L}\gpv\}$ 
is an orthonormal family in ${\cal H}$.
If $w\ne 0$, there exists $L_{0}\in {\bf N}^{*}$
such that $a_{L_{0}}\ne 0$.
By replacing $w$ by $a_{L_{0}}^{-1}s_{L_{0}}^{*}w$,
we can always assume $\langle w|\gpv\rangle=1$.
On the other hand,
$s_{J}^{*}s_{K}\gpv=0$ when $K\ne J^{n}$ for some $n\geq 0$.
Hence
$(s_{J}^{*})^{n}w\to \gpv$ when $n\to \infty$.
From this, we can obtain $\gpv$ from a given $w\in {\cal H}$,
$w\ne 0$.
This implies that ${\cal H}$ is irreducible.

Assume that $J$ is not nonperiodic (=periodic).
Then there exists $J_{0}\in {\bf N}^{*}$ such that 
$J$ is the $n$-times concatenation of $J_{0}$ for some $n\geq 2$.
Assume that a unit vector $\gpv$ satisfies $s_{J}\gpv=\gpv$.
Define two vectors $w_{1}$ and $w_{2}$ by
%
%
\begin{eqnarray}
\label{eqn:atwo}
w_{1}\equiv& \gpv+s_{J_{0}}\gpv+\cdots 
+s_{J_{0}}^{n-1}\gpv,\\
\nonumber
\\
\label{eqn:athree}
w_{2}\equiv& \gpv+\zeta s_{J_{0}}\gpv+\cdots 
+\zeta^{n-1} s_{J_{0}}^{n-1}\gpv
\end{eqnarray}
where $\zeta\equiv e^{2\pi\sqrt{-1}/n}$.
Since $s_{J_{0}}w_{1}=w_{1}$ and $s_{J_{0}}w_{2}=\zeta^{-1} w_{2}$,
$\langle w_{1}|w_{2}\rangle =0$.
Define
$V_{i}\equiv \overline{\coni w_{i}}$ for $i=1,2$.
Then we can verify that $V_{1}$ and $V_{2}$ are orthogonal.
Since  $\{0\}\ne V_{i}\subsetneq V_{1}\oplus V_{2}\subset {\cal H}$ for $i=1,2$,
${\cal H}$ is not irreducible.

\noindent
(b)
Assume $J\in {\bf N}^{\infty}$.
If $J$ is nonperiodic, then we can prove the irreducibility of $P(J)$ 
by the analogy of the case of (a).

Assume that $J$ is not nonperiodic.
We can assume that there exists $J_{0}\in {\bf N}^{*}$
such that $J$ is purely periodic with the repeating block $J_{0}$.
Assume that $\coni$ acts on a Hilbert space ${\cal K}$
with a cyclic unit vector $\gpv_{0}$ such that $s_{J_{0}}\gpv_{0}=\gpv_{0}$.
Let ${\cal L}$ denote the Hilbert space of 
all ${\cal K}$-valued functions $\phi $ on $U(1)\equiv\{z\in {\bf C}:|z|=1\}$
such that $\int_{U(1)}\|\phi(z)\|^{2}\,d\mu(z)<\infty$
where $\mu$ denotes the probabilistic Haar measure of the unitary group $U(1)$.
Define the action of $\coni$ on ${\cal L}$ by
%
%
\begin{equation}
\label{eqn:afour}
(s_{i}\phi)(z)\equiv zs_{i}\phi(z)
\quad(z\in U(1),\,\phi \in {\cal L},\,i\in {\bf N}).
\end{equation}
Let $\tilde{\gpv}_{0}$ denote the constant function in ${\cal L}$
such that $\tilde{\gpv}_{0}(z)\equiv \gpv_{0}$ for each $z\in U(1)$.
Then $(s_{J_{0}}\tilde{\gpv}_{0})(z)=z\tilde{\gpv}_{0}(z)$ for each $z\in U(1)$.
From this,
we see that $\{(s_{J_{(n)}})^{*}\tilde{\gpv}_{0}:n\geq 1\}$
is an orthonormal family in ${\cal L}$.
Let ${\cal L}_{1}$ denote the cyclic subspace of ${\cal L}_{1}$ generated by 
$\tilde{\gpv}_{0}$ with respect to the action of $\coni$.
Then ${\cal L}_{1}$ is $P(J)$.
Define the operator $T$ on ${\cal L}$ by
$(T\phi)(z)\equiv z\phi(z)$ for $\phi \in {\cal L}$ and $z\in U(1)$.
Then $T{\cal L}_{1}\subset {\cal L}_{1}$ and 
$Ts_{i}=s_{i}T$  and
$Ts_{i}^{*}=s_{i}^{*}T$ 
for each $i\in {\bf N}$ on ${\cal L}_{1}$.
Hence $T$ is the nontrivial element of the commutant of 
$\coni$ on ${\cal L}_{1}$.
Therefore neither ${\cal L}_{1}$ nor ${\cal L}$ is irreducible.
Hence $P(J)$ is not irreducible.
\qedh

%
%
\ssft{Proof of Proposition \ref{prop:shift}}
\label{subsection:apponetwo}

\noindent
(i)
For $[a]\in {\bf N}^{\infty}/\!\!\sim$,
we see that $\beta_{i}([a])\subset [a]$ and
$\beta_{i}^{-1}([a])\subset [a]$ for each $i$.
This implies that 
$\pi_{\beta}(\coni){\cal K}_{[a]} \subset {\cal K}_{[a]}$.
By definition, $l_{2}({\bf N}^{\infty})$ is decomposed into the direct sum of 
the family $\{{\cal K}_{[a]}:[a]\in {\bf N}^{\infty}/\!\!\sim\}$ of closed subspaces
as a Hilbert space, the statement of decomposition holds. 
It is sufficient to show the irreducibility of ${\cal K}_{[a]}$
for each $[a]$.
We consider the following two cases.

\noindent
(a) Assume that $a$ has the repeating block $a^{'}$.
Let $(a^{'})^{\infty}$ denote the purely periodic sequence with 
the repeating block $a^{'}$.
Then
$[a]=\{\beta_{J}(a^{'}):J\in {\bf N}^{k},\,k\geq 0\}$
and $\beta_{a^{'}}((a^{'})^{\infty})=(a^{'})^{\infty}$.
This implies that $({\cal K}_{[a]},\pi_{\beta}|_{{\cal K}_{[a]}})$ is $P(a^{'})$
with the GP vector $e_{(a^{'})^{\infty}}$.
From Proposition \ref{prop:oniprop}(v)
and the nonperiodicity of $a^{'}$,
${\cal K}_{[a]}$ is irreducible.

\noindent
(b) Assume that $a$ has no repeating block.
By the analogy of the case of (a),
we see that  ${\cal K}_{[a]}$ is $P(a)$ with the GP vector $e_{a}$.
From Proposition \ref{prop:oniprop}(v)
and the nonperiodicity of $a$,
${\cal K}_{[a]}$ is irreducible.

\noindent
(ii)
From the proof of (i),
$\theta_{[a]}$ is $P(J)$ for a nonperiodic element 
$J$ in ${\bf N}^{*}\cup {\bf N}^{\infty}$.
From this and Proposition \ref{prop:oniprop}(iii) and (iv),
the statement holds.

\noindent
(iii)
The proof has been already given in the proof of (i).

\noindent
(iv)
Let $({\cal H},\pi)$ be an irreducible permutative representation of $\coni$.
From Proposition \ref{prop:oniprop}(ii),
$({\cal H},\pi)$ is $P(J)$
for some nonperiodic element $J$ in ${\bf N}^{*}\cup {\bf N}^{\infty}$.
From (iii),
if $J\in {\bf N}^{*}$,
then 
$\theta_{[a]}$ is $P(J)$ for $a\equiv J^{\infty}\in {\bf N}^{\infty}$, 
and if $J\in {\bf N}^{\infty}$,
then $\theta_{[J]}$ is $P(J)$.
Hence the statement holds.
\qedh


%
\end{document}